\documentstyle[12pt]{article}
\topmargin  - 0.1cm
\newcommand\R{\hbox{I\kern - .2em\hbox{R}}}
\newcommand\N{\hbox{I\kern - .2em\hbox{N}}}

\newtheorem{dfn}{Definition}[section]

\newtheorem{theorem}[dfn]{Theorem}
\newtheorem{cor}[dfn]{Corollary}
\newtheorem{lem}[dfn]{Lemma}
\def\blacksquare{\hbox{\vrule width 4pt height 4pt depth 0pt}}
\def\abstract{\if@twocolmn\section*{Abstract}\else \small \begin{center}
{\bf Abstract\vspace{ - .5em}\vspace{0pt}}
\end{center}
\quotation
\fi}

\begin{document}

\title{{\LARGE {$\left(\varphi_1,  \varphi_2\right)-$Variational principle}}}

\author{Abdelhakim MAADEN and Abdelkader STOUTI}

\date{{Universit\'e Sultan Moulay Slimane \\  Facult\'e des Sciences et Techniques.  Laboratoire de  Math\'ematiques et Applications.  B.P. 523. BENI - MELLAL,  MAROC.}\\ {Corresponding author: E-mail address:  hmaaden2002@yahoo.fr}}

\maketitle

\baselineskip=20pt

\begin{abstract}
In this paper we prove that if $X $ is a Banach space, then for every lower semi-continuous bounded below function $f, $ there exists a $\left(\varphi_1, \varphi_2\right)-$convex function $g, $ with arbitrarily small norm,  such that $f + g $ attains its strong minimum on $X. $ 
This result extends some of the  well-known varitional principles as that of Ekeland [18], that of Borwein-Preiss [6] and that of Deville-Godefroy-Zizler [14, 15].
\end{abstract}



\section{Introduction}
Let $\left(X, \| . \|\right) $ be a Banach space. Let $f $ be a real-valued function defined on $X, $ lower semicontinuous and bounded below. Let $P $ be a class of functions in $X $ which serves as a source of possible perturbations for $f. $ By a variational principle we mean an assertion ensuring the existence of at least one perturbation $g $ from $P $ such that $f + g $ attains its minimum on $X. $

The first variational principle, based on the Bishop-Phelps lemma [3, 27] , was established by Ekeland  [18]. It says that $P $ is a family of suitable translations of the norm. 

If $g$ is required to be smooth, then we speak about the smooth variational principle. The first result of this type was shown by Stegall [31, 27], where $P $ is the elements of the dual space $X^*. $ He proved that if $X $ has the Radon-Nikodym property in particular, if $X $ is reflexive, then one can take for  $g $ even a linear functional, with arbitrarily small norm. In [17], Deville-Maaden shown under the same hypotheses as Stegall that $X $ has the Radon-Nikodym property and  the function $f $ is lower semicontinuous and super-linear with the set $P$ consisting of radial smooth functions.
However, this principle does not cover some important Banach spaces. For example the space $c_0$ does not have the Radon-Nikodym property while it, in fact, admits a smooth norm [5]. In this direction Borwein-Preiss [6] proved a smooth variational principle imposing no additional conditions on the space, except, the presence of some equivalent smooth norm with the set $P $ being the family of smooth combinations of the norm in $X.$ Haydon's work [23], shows that there exists a Banach spaces with smooth bump function without an equivalent smooth norm (a function $b $ is bump if it has a non empty and  bounded support). So,  the variational principle of Borwein-Preiss is not applicable in these spaces. So that, Deville et al [14, 15] extended the Borwein-Preiss variational principle to spaces with smooth bump function, with $P $ equal to the family of Lipschitz smooth functions.

In an analytical approach we can often associate a geometrical approach to complete study of which or stimulates the analytical approche. From this perspective Browder [8] gaves a geometrical result  which bears at present the name of the Drop Theorem  (see also [10]). Penot in [26, 21] showed that the drop theorem is a geometrical version of the Ekeland's variational principle. After this, Maaden in [25, 22] introduced and studied the notion of the smooth drop which can be seen as a geometrical version of the smooth variational principle of Borwein-Preiss.

Those variational principles are a tools that have been very important in nonlinear analysis, in that they enjoyed a big deal of applications from the geometry of Banach spaces [3, 4, 7] to the optimization theory [18, 19, 30]  and of generalized differential and sub-differential calculus [1, 2, 6, 11, 12, 13, 26], calculus of variations [9, 18] up to the nonlinear semi-groups theory [7, 18] and the viscosity solutions of Hamilton-Jacobi equations [12, 13, 15].

In [28, 29], Pini et all defined the notion of $\left(\varphi_1, \varphi_2\right)$-convex functions. They say that a real valued function $f $ defined on a non empty subset $D $ of $\R^n $ is 
$\left(\varphi_1, \varphi_2\right)-$convex if 
$f\left(\varphi_1\left(x, y, \lambda\right)\right) \le \varphi_2\left(x, y, \lambda, f\right) $ for all $x, y \in D $ and for all $\lambda \in [0, 1].$ Where $\varphi_1 $ is a function from $D \times D \times [0, 1] $ in $\R^n $ and $\varphi_2 $ is a function from $D \times D \times [0, 1] \times F$ in $\R, $ with $F $ is a given vector space of real valued functions defined on the set $D. $ In this paper we shall use the same definition of $\left(\varphi_1, \varphi_2\right)$-convex functions  as above with using any Banach spaces instead of $\R^n.$  In this way, we prove that under suitable choices of the functions $\varphi_1$ and $\varphi_2$ a new variational principle for the set of $\left(\varphi_1, \varphi_2\right)$-convex functions (see Theorem 3.1). This $\left(\varphi_1, \varphi_2\right)-$ variational principle is providing a unified framework to deduce Ekeland's, Borwein-Preiss's and Deville's variational principles.

\section{Auxiliaries results}

In this section we shall give some definitions and establish some auxiliaries results which we shall use to prove our main result in this paper.

Let $\left(X, \left\|.\right\|\right) $ be a Banach space. For a continuous function $f: X \longrightarrow \R $ we define
$$\mu \left(f\right) = \sum_{n = 1}^\infty \displaystyle\frac{\left\|f\right\|_n}{2^n},  $$
where
$$\left\|f\right\|_n = \sup\left\{|f\left(x\right)| ; x \in X, \left\|x\right\| \le n\right\}. $$
Let $M $ be the set of all continuous functions $f $ such that $\mu\left(f\right) < \infty. $ It is routine to check that $\left(M, \mu\right) $ is a Banach space.

Let $\varphi_1: X \times X \times [0, 1]  \longrightarrow X \hbox{ and } \varphi_2: X \times X \times [0, 1] \times F \longrightarrow \R, $ two functions where $F $ is a given set of real functions on $X. $
Define,

\begin{dfn} A function $f: X  \longrightarrow \R $ is said to be  $\left(\varphi_1, \varphi_2\right)-$convex  if
$$f\left(\varphi_1\left(x, y, \lambda\right)\right) \le \varphi_2\left(x, y, \lambda, f\right), \forall x, y \in X, \forall \lambda \in [0, 1].$$
\end{dfn}

Remarking that under suitable assumptions on $\varphi_1$ and/or $\varphi_2, $ the class of $\left(\varphi_1, \varphi_2\right)-$ convex functions is a convex cone. For example:

\noindent 1) If $\varphi_2 $ is super-linear with respect to $f \in F $ (that $\varphi_2 $ is super-additive and positively homogeneous), then the class of $\left(\varphi_1, \varphi_2\right)-$convex functions is a convex cone.

Indeed, let $f, g $ are two $\left(\varphi_1, \varphi_2\right)-$convex functions and $\alpha > 0. $ Then, for $x, y \in X $ and $\lambda \in [0, 1] $ we have
\begin{eqnarray*}\left(f + g\right)\left(\varphi_1\left(x, y, \lambda\right)\right)  &\le &  \varphi_2\left(x, y, \lambda, f\right) + \varphi_2\left(x, y, \lambda, g\right) \\
&\le & \varphi_2\left(x, y, \lambda, f + g\right)
\end{eqnarray*}
and
\begin{eqnarray*}\left(\alpha f\right)\left(\varphi_1\left(x, y, \lambda\right)\right)  &= &  \alpha\left(f\left(\varphi_1\left(x, y, \lambda\right)\right)\right)  \\
&\le & \alpha\varphi_2\left(x, y, \lambda, f\right)\\
&= & \varphi_2\left(x, y, \lambda, \alpha f\right).
\end{eqnarray*}

\noindent 2) If $\varphi_2\left(x, y, \lambda, f\right) = C\left(\left(1 - \lambda\right)f\left(x\right) + \lambda f\left(y\right)\right) $ for some  $C > 0, $ the set of $\left(\varphi_1, \varphi_2\right)-$convex functions is a convex cone.

\bigskip

In all the  sequel, we define the following sets :
$$\Phi = \left\{f \in M : f \hbox{ is } \left(\varphi_1, \varphi_2\right)-\hbox{convex and } f \ge 0\right\},$$
$$F = \left\{f \in \Phi : f\left(x\right) \longrightarrow +\infty \hbox{ as } ||x|| \longrightarrow +\infty\right\}. $$
The metric $\rho $  on $\Phi $ is defined as:
$$\rho\left(f, g\right) = \mu\left(f - g\right)= \sum_{n\ge1}\frac{\|f - g\|_n}{2^n} \hbox { for all } f, g \in \Phi, $$
and it is easy to show that $\left(\Phi, \rho\right) $ is a complete metric space.

Throughout this paper, the functions $\varphi_1 $ and $\varphi_2 $     satisfies the following assumptions:

\noindent $\left(P_1\right)$ $\varphi_1\left(x, x, 0\right) = x ; \forall x \in X $ ; 

\noindent $\left(P_2\right)$ $\varphi_1\left(x, y, \lambda\right) + \varphi_1\left(z, z, 0\right) = \varphi_1\left(x + z, y + z, \lambda\right) ; \forall x, y, z  \in X,   \forall \lambda \in [0, 1] $ ;

\noindent $\left(P_3\right)$ $\exists C \ge 1, $ such that $\varphi_2\left(\lambda x, \lambda x, 0, h\right) \le C[\left(1 - \lambda\right)h\left(0\right) + \lambda h\left(x\right)] ; \forall x \in X, \forall \lambda \in [0, 1], \forall h \in \Phi $  ;

\noindent $\left(P_4\right)$ For $x_0 \in X, $ $\varphi_2\left(x - x_0, y - x_0, \lambda, h\right) \le \varphi_2\left(x, y, \lambda, h\left(. - x_0\right)\right) ; \forall x, y \in X, \forall \lambda \in [0, 1] ; \forall h \in \Phi $  ;

\noindent $\left(P_5\right)$ The class of $\left(\varphi_1, \varphi_2\right)-$ convex functions is a convex cone.

\noindent  We will also assume that $\varphi_1 $ is continuous with respect to $\lambda. $

\medskip

We present now three preliminaries lemmas, which are useful for the proof of our principal result of this paper. In the first, we use $\left(P_1\right) $ and  $\left(P_3\right) $  to prove the following:

\begin{lem}
Let $h \in \Phi $   and let $y = \lambda x, \, \lambda > 1. $ Then, $h\left(y\right) - h\left(0\right) \ge \displaystyle\frac{\lambda}{C}\left(h\left(x\right) - Ch\left(0\right)\right). $
\end{lem}

\noindent Proof:

Let $\mu = 1/\lambda. $ Then $x = \mu y. $ By using $\left(P_1\right) $ and  $\left(P_3\right) $  we obtain that:
\begin{eqnarray*}
h\left(x\right) &= & h\left(\mu y\right) \\
&= &  h\left(\varphi_1\left(\mu y, \mu y, 0\right)\right)  \\
&\le & \varphi_2\left(\mu y, \mu y, 0, h\right) \\
&\le & C\left(\left(1 - \mu\right)h\left(0\right) + \mu h\left(y\right)\right).
\end{eqnarray*}
Consequently, we get 
$$h\left(x\right) - Ch\left(0\right) \le C\mu \left(h\left(y\right) - h\left(0\right)\right). $$
Since $c >0 $ and $\mu > 0, $ we deduce that
$$h\left(y\right) - h\left(0\right) \ge \displaystyle\frac{1}{C\mu}\left(h\left(x\right) - Ch\left(0\right)\right) = \displaystyle\frac{\lambda}{C}\left(h\left(x\right) - Ch\left(0\right)\right) $$
and the proof is complete.
\blacksquare

We have now all tools to confirm that $\left(F, \rho\right) $ is a Baire space. For this, it suffices  to show that

\begin{lem}
$\left(F, \rho\right) $ is open in $\Phi. $
\end{lem}

\noindent Proof:

Let $f $ in $F. $ Let $N > Cf\left(0\right) $ (where $C $ is given by $\left(P_3\right)$) and let $\varepsilon > 0 $ be such that
$$0 < \varepsilon < \displaystyle\frac{N - Cf\left(0\right)}{2C + 1}. \eqno{\left(1\right)}$$
Since $f \in F ; f \longrightarrow +\infty $ as $\left\|x\right\| \longrightarrow +\infty ; $ there exists $n \in \N $ such that
$$f\left(x\right) > N \hbox{ whenever  } \left\|x\right\| \ge n. \eqno{\left(2\right)} $$
Let $g \in \Phi $ such that $\rho\left(f, g\right) < \displaystyle\frac\varepsilon {2^n}. $
Then
$$\displaystyle\frac{\|f - g\|_n}{2^n} < \displaystyle\frac\varepsilon {2^n}. $$
Hence,
$$|f\left(x\right) - g\left(x\right)| < \varepsilon \quad\hbox{ whenever}\quad \left\|x\right\| \le n. \eqno{\left(3\right)} $$
In particular, we have
$$|f\left(0\right) - g\left(0\right)| < \varepsilon. \eqno{\left(4\right)} $$
Combining (2) and (3) we obtain for $\left\|x\right\| = n, $
$$g\left(x\right)  >  N - \varepsilon. \eqno{\left(5\right)} $$
On the first, for $y \in X $ such that $\left\|y\right\| \longrightarrow +\infty, $ there exist $x \in X $ with $\left\|x\right\| = n $ and $\lambda > 1 $ such that $y = \lambda x, $ and we have $\lambda \longrightarrow +\infty. $
 Therefore, combining Lemma 2.2, (5), (4) and (1)  we obtain
\begin{eqnarray*}g\left(y\right) - g\left(0\right) &\ge & \displaystyle\frac{\lambda}{C}\left(g\left(x\right) - Cg\left(0\right)\right) \\
&> & \displaystyle\frac{\lambda}{C}\left(N - \varepsilon - Cg\left(0\right)\right) \\
&\ge & \displaystyle\frac{\lambda}{C}\left(N - \varepsilon - Cf\left(0\right) - C\varepsilon\right) \\
&> & \lambda \varepsilon,
\end{eqnarray*}
which implies that
$$g\left(y\right) > \lambda\varepsilon + g\left(0\right) \ge \lambda\varepsilon > 0. $$
On the other hand,  we know that $\lambda \longrightarrow +\infty $ as $\left\|y\right\| \longrightarrow +\infty. $  So, we get that $g\left(y\right) \longrightarrow +\infty $ as $\left\|y\right\| \longrightarrow +\infty. $ So that,  $g \in F $ and $F $ is open.
\blacksquare

Next, by using $\left(P_1\right), \left(P_2\right) $ and $\left(P_4\right) $ we obtain the following:

\begin{lem}
Let $\theta $ be a $\left(\varphi_1, \varphi_2\right)-$convex function and let $h\left(x\right) = \theta\left(x - x_0\right). $ Then, $h $ is a $\left(\varphi_1, \varphi_2\right)-$convex function.
\end{lem}

\noindent Proof:

Let $x, y \in X $ and $\lambda \in [0, 1]. $ By using $\left(P_1\right), \left(P_2\right) $ and $\left(P_4\right)  $ we get
\begin{eqnarray*}
h\left(\varphi_1\left(x, y, \lambda\right)\right) &= & \theta\left(\varphi_1\left(x, y, \lambda\right) - x_0\right) \\
&= & \theta\left(\varphi_1\left(x , y, \lambda\right) + \varphi_1\left(-x_0, -x_0, 0\right)\right) \\
&= & \theta\left(\varphi_1\left(x - x_0, y - x_0, \lambda\right)\right) \\
&\le & \varphi_2\left(x - x_0, y - x_0, \lambda, \theta\right) \\
&\le & \varphi_2\left(x, y, \lambda, \theta\left(. - x_0\right)\right) \\
&= & \varphi_2\left(x, y, \lambda, h\right),
\end{eqnarray*}
which shows that $h $ is a $\left(\varphi_1, \varphi_2\right)-$convex function.
\blacksquare

\begin{cor}
Let $\theta $ be a $\left(\varphi_1, \varphi_2\right)-$convex function in $F $ then the function  $h\left(x\right) = \theta\left(x - x_0\right)$  is in $F. $
\end{cor}

\section{The main result }
In this section we shall establish a $(\varphi_1, \varphi_2)-$variational principle. We show that the set $P $ ; whish is a source of perturbation for $f $ ;  is a class of $(\varphi_1, \varphi_2)-$convex functions. Furthermore we can take them of  $C^\infty $ in smooth Banach spaces.

In the mathematical field of topology, a $G_\delta $ set is a subset of a topological space that is a countable intersection of open sets. 
In a complete metric space, a countable union of nowhere dense sets is said to be meagre; the complement of such a set is a residual set. 

An element $y$ of a Banach space $X $ is said a strong minimum for a real function $f $ defined on the space $X, $ if $f(y) $ is the infimum of $f $ and any minimizing sequence for $f $ converges to $y. $

The aim result in this paper is the following variational principle:

\begin{theorem}
Let $X $ be a Banach space. Let $f: X \longrightarrow \R \cup \left\{+\infty\right\} $ be a lower semi-continuous function bounded from below. Let $Y $ be a subset of $F $ such that:

i) the metric $\rho_Y $ in $Y $ is such that $\rho_Y\left(f, g\right) = \mu_Y\left(f - g\right) \ge \mu \left(f - g\right), $ for all $f, g \in Y.$

ii) $\left(Y, \rho_Y\right) $ is a Baire space.

iii)  there exists $ \theta \in Y $  such that $\mu_Y\left(\theta\right) < +\infty, \theta\left(0\right) = 0, $ there is $k \in ]0, 1[ $ such that  for every  $ \left\|x\right\| \ge k $ we have $\theta\left(x\right) \ge k^2 $ and $\mu_Y\left(\theta\left(. - x_0\right)\right) \le \mu_Y\left(\theta\right) + ||\theta||_{||x_0||}. $

Then the set
$$\left\{g \in Y : f + g \hbox{ attains its strong minimum on } X\right\}
$$
is residual in $Y. $
\end{theorem}

 Next, we shall show that Theorem 3.1 is providing a unified framework to deduce Ekeland's variational principle [18], Borwein-Preiss's [6]  variational principle and Deville-Godefroy-Zizler's Varitional principle [15].

\noindent{\bf Application 1.}
As a first application we get the Ekeland's variational principle [18].

Let $\left(X, \left\|.\right\|\right) $ be a Banach space. Assume that
$\varphi_1\left(x, y, \lambda\right) = \lambda x + \left(1 - \lambda\right)y $ and $\varphi_2\left(x, y, \lambda, f\right) = \lambda f\left(x\right) + \left(1 - \lambda\right)f\left(y\right). $ Then $\varphi_1 $ and $\varphi_2 $ satisfies $\left(P_1\right), \left(P_2\right), \left(P_3\right)$ and $\left(P_4\right). $
Let
$$Y = \left\{f: X \longrightarrow \R : f  \hbox { convex, Lipschitz, } \ge 0, f \longrightarrow + \infty \hbox { as } ||x|| \longrightarrow + \infty\right\}. $$
We define on $Y $ the metric $\rho_Y $ such that for $f, g \in Y, $
$$\rho_Y\left(f, g\right) = \mu_Y\left(f - g\right) = \sum_{n \ge 1}\displaystyle\frac{||f - g||_n}{2^n} + \sup\left\{\displaystyle\frac{|\left(f - g\right)\left(x\right) - \left(f - g\right)\left(y\right)|}{||x - y||} ; x \not= y\right\}. $$
It is clear that $\left(Y, \rho_Y\right) $ satisfies $\left(P_5\right)$ and the conditions $(i) $ and $ (ii)$ of Theorem 3.1. Also, the function $\theta = ||x|| $ satisfies the assertion $(iii) $ of Theorem 3.1. Consequently we have the following:

\begin{cor}
Let $\left(X, ||.||\right) $ be a Banach space,  consider a lower semi-continuous bounded below function $f: X \longrightarrow \R \cup \left\{+\infty\right\}. $ Then for each $\varepsilon > 0, $ there exists $x_0 \in X $ such that $$f\left(x\right) + \varepsilon||x - x_0|| \ge f\left(x_0\right). $$
\end{cor}

\noindent Proof:

From Theorem 3.1, for each $\varepsilon > 0, $ there exits $g \in Y $ such that $\mu_Y\left(g\right) < \varepsilon $ and $f + g $ attains a strong minimum at $x_0. $ Therefore, for all $x \in X, $
$$f\left(x\right) + g\left(x\right) \ge f\left(x_0\right) + g\left(x_0\right)
\hbox{ and }
\sum_{n \ge 1}\displaystyle\frac{||g||_n}{2^n} + \sup\left\{\displaystyle\frac{|g\left(x\right) - g\left(y\right)|}{||x - y||} ; x \not= y\right\} < \varepsilon, $$
which implies that
\begin{eqnarray*}
f\left(x\right) &\ge & f\left(x_0\right) + g\left(x_0\right) - g\left(x\right) \\
&\ge & f\left(x_0\right) - \varepsilon||x - x_0||.
\end{eqnarray*}
\blacksquare

\noindent{\bf Application 2.} Let $\left(X, \left\|.\right\|\right) $ be a Banach space with smooth norm.
Assume that $\varphi_1\left(x, y, \lambda\right) = \lambda x + \left(1 - \lambda\right)y $ and $\varphi_2\left(x, y, \lambda, f\right) = \lambda f\left(x\right) + \left(1 - \lambda\right)f\left(y\right). $ Then $\varphi_1 $ and $\varphi_2 $ satisfies $\left(P_1\right), $$\left(P_2\right), $$\left(P_3\right), $ and $\left(P_4\right). $
Let
$$Y = \left\{f \hbox{ is } C^1\hbox{-smooth, convex,} \ge 0 \hbox{ and } f \longrightarrow +\infty \hbox{ as } ||x|| \longrightarrow +\infty\right\}. $$
We define the metric $\rho_Y $ in $Y $ by:
$$\rho_Y(f, g) = \mu_Y (f - g) = \sum_{n \ge 1}\displaystyle\frac{||f - g||_n}{2^n} + ||(f - g)'||_\infty \hbox{ for all } f, g \in Y$$
where $\left\|f'\right\|_{\infty} := \displaystyle\sup_{\|x\| \le 1}\left\|f'(x)\right\|_{X^*} $  and the space $\left(Y, \rho_Y\right) $ satisfies $(i) $ and $(ii) $ of  Theorem 3.1 and so also $\left(P_5\right). $

Let
\begin{eqnarray*}
h : &[0, +\infty[ & \longrightarrow [0, +\infty[ \\
&t & \longmapsto \cases{t^2 \quad \quad\phantom{rr} \hbox{if }\quad 0 \le t \le 1 \cr
2t - 1 \quad \hbox{if} \quad \phantom{r}t > 1.}
\end{eqnarray*}
The function $\theta\left(x\right) = h\left(||x||\right) \in Y $ satisfies the assertion $\left(iii\right) $ of Theorem 3.1.

Therefore we have the Borwein-Preiss's variational principle [6, 27]:

\begin{cor}
Let $\left(X, ||.||\right) $ be a Banach space with a smooth norm and consider a lower semi-continuous function $f: X \longrightarrow \R \cup \left\{+\infty\right\}$ bounded from below.  Then the set
$$\left\{g \in Y : f + g \hbox{ attains its strong minimum on } X\right\}$$
is residual in $Y. $
\end{cor}

\noindent{\bf Application 3.}
Let $X $ be a Banach space admitting Lipschitz $C^1-$smooth bump function.
According to a construction of Leduc [24], there exists a Lipschitz  function $d: X \longrightarrow \R $ which is $C^1-$smooth on $X \setminus \left\{0\right\} $  and satisfies:

i) $d\left(\lambda x\right) = \lambda d\left(x\right) $ for all $\lambda > 0 $ and for all $x \in X, $

ii) there exists $C \ge 1 $ such that
$\left\|x\right\| \le d\left(x\right) \le C \left\|x\right\| \hbox{ for all } x \in X. $

\noindent Moreover the function $d^2 $ is $C^1-$smooth on all the space $X. $

Let $\varphi_1\left(x, y, \lambda\right) = \lambda x + \left(1 - \lambda\right)y $ and $\varphi_2\left(x, y, \lambda, f\right) = C^2[\lambda f\left(x\right) + \left(1 - \lambda\right)f\left(y\right)]. $ Then $\varphi_1 $ and $\varphi_2 $ satisfies
$\left(P_1\right), \left(P_2\right), \left(P_3\right)$ and $\left(P_4\right). $
Let
$\theta\left(x\right) = d^2\left(x\right). $ We have
$$d^2\left(\lambda x + \left(1 - \lambda\right)y\right) \le C^2||\lambda x + \left(1 - \lambda\right)y||^2. $$
Since the function $||.||^2 $ is convex we deduce that
$$d^2\left(\lambda x + \left(1 - \lambda\right)y\right) \le C^2\left(\lambda d^2\left(x\right) + \left(1 - \lambda\right)d^2\left(y\right)\right).$$
That is the function $d^2 $ is a $\left(\varphi_1, \varphi_2\right)-$convex function.

Let
$$Y = \left\{f \hbox{ a } \left(\varphi_1, \varphi_2\right)-\hbox{convex}, C^1-\hbox{Lipschitz,}  \ge 0 \hbox{ and } f \longrightarrow +\infty \hbox{ as } ||x|| \longrightarrow +\infty\right\}$$
and so the set $Y $ satisfies $\left(P_5\right). $

The metric $\rho_Y $ on $Y $ is such that, for $f, g \in Y$
$$\rho_Y(f, g) = \mu_Y(f - g) = \sum_{n \ge 1}\displaystyle\frac{||f - g||_n}{2^n} + \displaystyle\sum_{n \ge 1}\frac{||(f - g)'||_n}{2^n} $$ where 
$\left\|f^{'}\right\|_n = \displaystyle\sup_{||x|| \le n}\left\|f^{'}(x)\right\|_{X^{*}}. $

In the other hand,
let $\theta\left(x\right) = d^2\left(x\right). $ So that,

i) $\theta\left(0\right) = 0$

ii) $\mu_Y\left(\theta\right) < \infty $

iii) let $0 < k < 1. $ Hence, for all $x \in X $ such that $||x|| \ge k $ we have $d^2\left(x\right) \ge ||x||^2 \ge k^2. $

\noindent Therefore the function $\theta \in Y $ and  satisfies $\left(iii\right) $ of Theorem 3.1.

Thus we have the following variational principle (for unbounded functions) of Deville-Godefroy-Zizler  [14, 15, 16, 20]:

\begin{cor}
Let $\left(X, ||.||\right) $ be a Banach space admitting a $C^1-$Lipschitz bump function and consider a lower semi-continuous bounded below function $f: X \longrightarrow \R \cup \left\{+\infty\right\}. $ Then the set
$$\left\{g \in Y : f + g \hbox{ attains its strong minimum on } X\right\}$$
is residual in $Y. $
\end{cor}

Now, we are ready to give the proof of Theorem 3.1.

\begin{center}
{\bf Proof of Theorem 3.1}
\end{center}

Following the method of [15, 20], for $n \in \N \setminus \left\{0\right\}, $ we let
$$G_n = \left\{g \in Y : \exists x_0 \in X, \left(f + g\right)\left(x_0\right) < \inf\left\{\left(f + g\right)\left(x\right) : \left\|x - x_0\right\| \ge 1/n\right\}\right\}. $$

\noindent Claim 1. We claim that $G_n $ is open for each $n. $ Indeed,
let $n \in \N $ and  $g \in G_n. $ So that there is $x_0 $ in $X $ such that
$$\left(f + g\right)\left(x_0\right) < \inf\left\{\left(f + g\right)\left(x\right) : \left\|x - x_0\right\| \ge 1/n\right\}.$$
Let $0 < \varepsilon <1 $ such that
$$\left(f + g\right)\left(x_0\right) + 2\varepsilon < \inf\left\{\left(f + g\right)\left(x\right) : \left\|x - x_0\right\| \ge 1/n\right\}. \eqno{\left(1\right)} $$
Let $A = C\left(f + g\right)\left(x_0\right) + C\left(g\left(0\right) - \inf\left(f\right)\right) + \left(2C + 3\right)\varepsilon, $ where $C $ is given by $\left(P_3\right). $
Since $g \in Y, $
$g $ goes to $+\infty $ as $\left\|x\right\| $ goes to $+\infty. $ This means that,
there is  $k $ in $\N $ such that $k > \left\|x_0\right\| $ and $g\left(x\right) > A $ whenever $\left\|x\right\| \ge k. $
This is equivalent to say that
$$g\left(x\right)  > C\left(f + g\right)\left(x_0\right) + C\left(g\left(0\right) - \inf\left(f\right)\right) + \left(2C + 3\right)\varepsilon  \quad \hbox{ whenever } \left\|x\right\| \ge k. \eqno{\left(2\right)} $$
Let $h \in Y $ such that $\rho_Y\left(h, g\right) < \displaystyle\frac{\varepsilon}{2^k}. $ We have
$$\sum_{n\ge 1}\displaystyle\frac{\left\|h - g\right\|_n}{2^n} \le \rho_Y\left(h, g\right) = \mu_Y \left(h - g\right)  < \displaystyle\frac{\varepsilon}{2^k}.  $$
Thus
$$\displaystyle\frac{\|h - g\|_k}{2^k} < \displaystyle\frac{\varepsilon}{2^k}. $$
So  that
$$|h\left(x\right) - g\left(x\right)| < \varepsilon \quad \hbox{ whenever } \left\|x\right\| \le k, \eqno{\left(3\right)} $$
in particular
$$|h\left(x_0\right) - g\left(x_0\right)| < \varepsilon.  \eqno{\left(4\right)} $$
Combining (2) with (3) we obtain that
$$h\left(x\right) > C\left(f + g\right)\left(x_0\right) + C\left(g\left(0\right) - \inf\left(f\right)\right) + \left(2C + 2\right)\varepsilon > 0 \quad \hbox{ whenever } \left\|x\right\| = k. $$
Since $C \ge 1 $ and $h \ge 0, $ we deduce that for $\|x\| = k $
$$h\left(x\right) \ge \displaystyle\frac{h\left(x\right)}{C} > \left(f + g\right)\left(x_0\right) + \left(g\left(0\right) - \inf\left(f\right)\right) + \left(2 + \left(2/C\right)\right)\varepsilon.   \eqno{\left(5\right)}$$

In the first hand,
let $y \in X $ such that $\left\|y\right\| > k. $ Then, there exist $\lambda > 1 $ and $x \in X $ with $\left\|x\right\| = k, $ such that $y = \lambda x. $
By using Lemma 2.2, we deduce that
$$h\left(y\right) - h\left(0\right) \ge  \displaystyle\frac{\lambda}{C}\left(h\left(x\right) - Ch\left(0\right)\right) \ge \displaystyle\frac{1}{C}\left(h\left(x\right) - Ch\left(0\right)\right) = \displaystyle\frac{h\left(x\right)}{C} - h\left(0\right). $$
Combining this with (5) we show for $\|y\| \ge k $ that,
$$h\left(y\right) - h\left(0\right) >  \left(f + g\right)\left(x_0\right) + g\left(0\right) - \inf f + \left(2 + \displaystyle\frac{2}{C}\right)\varepsilon - h\left(0\right). \eqno{\left(6\right)}$$
 Combining the fact that $h \ge 0, \left(6\right), \left(3\right) $ and $\left(4\right) $ we obtain for all $x \in X $ such that $\|x\| \ge k $:
\begin{eqnarray*}\left(f + h\right)\left(x\right)  &\ge &  \inf\left(f\right) + h\left(x\right) \\
&\ge & \inf\left(f\right) + h\left(x\right)- h\left(0\right)   \\
&> & \inf\left(f\right) + \left(f + g\right)\left(x_0\right) + g\left(0\right) - \inf\left(f\right)
+ \left(2 + \displaystyle\frac{2}{C}\right)\varepsilon -  h\left(0\right)
\\
&> & \left(f + g\right)\left(x_0\right) + \left(1 + \displaystyle\frac{2}{C}\right)\varepsilon
\\
&> & \left(f + h\right)\left(x_0\right) + \displaystyle\frac{2}{C}\varepsilon  \\
&> & \left(f + h\right)\left(x_0\right).
\end{eqnarray*}
Therefore for all $x \in X $ such that $\left\|x\right\| \ge k, $ we have
$$\left(f + h\right)\left(x\right) > \left(f + h\right)\left(x_0\right). $$

In other hand, if $\left\|x\right\| \le k $ and $\left\|x - x_0\right\| \ge 1/n, $ and combining (4), (1) and (3) we obtain that
\begin{eqnarray*}\left(f + h\right)\left(x_0\right)  &< & \left(f + g\right)\left(x_0\right) + \varepsilon \\
&\le & \inf\left\{\left(f + g\right)\left(x\right) : \left\|x - x_0\right\| \ge 1/n\right\} - 2\varepsilon + \varepsilon \\
&\le & \left(f + g\right)\left(x\right)  -  \varepsilon \\
&< & \left(f + h\right)\left(x\right).
\end{eqnarray*}
Then for all $x $ such that $\left\|x - x_0\right\| \ge 1/n $ we have
$$\left(f + h\right)\left(x_0\right) < \left(f + h\right)\left(x\right). $$
Hence $h \in G_n $ and $G_n $ is open.

Claim 2. We confirm that $G_n $ is dense in $Y. $ Indeed,
let  $g \in Y $ and $0 < \varepsilon < 1. $ Let $c > 0 $ be such that
$$\left(f + g\right)\left(x\right) > \inf\left(f + g\right) + 1 \hbox { whenever } \left\|x\right\| > c. $$
Let $1 > \delta > 0 $ be such that $\delta\left(\mu_Y\left(\theta\right) + \left\|\theta\right\|_c\right) < \varepsilon. $ Let $x_0 \in X $ be such that
$$\left(f + g\right)\left(x_0\right) < \inf\left(f + g\right) + \displaystyle\frac{\delta}{n^2}. \eqno{\left(6\right)} $$
Since $\displaystyle\frac{\delta}{n^2} < 1, $ we deduce that
$$\left\|x_0\right\| \le c. \eqno{\left(7\right)} $$
Let $h\left(x\right) = \delta\theta\left(x - x_0\right). $ Now Corollary 2.5 ensure that  $h $ is a $\left(\varphi_1, \varphi_2\right)-$convex  function in $F. $

From the hypothesis $(iii) $ of Theorem 3.1  and $(7), $ we get
$$\rho_Y\left(h, 0\right) = \mu_Y\left(h\right) = \delta\mu_Y\left(\theta\left(. - x_0\right)\right) \le \delta\mu_Y\left(\theta\right) + \delta\left\|\theta\right\|_{\left\|x_0\right\|} \le \delta\left(\mu_Y\left(\theta\right) + \left\|\theta\right\|_c\right) < \varepsilon. $$

Now if $\left\|x - x_0\right\| \ge 1/n, $ and by $\left(iii\right) $ of Theorem 3.1 we deduce that,
$$
h\left(x\right) =  \delta\theta\left(x - x_0\right) \ge  \displaystyle\frac{\delta}{n^2}.
$$
By using (6), we get
\begin{eqnarray*}
\inf\left\{f + g + h : \left\|x - x_0\right\| \ge 1/n\right\}&\ge &
\inf\left\{f + g  : \left\|x - x_0\right\| \ge 1/n\right\} + \displaystyle\frac{\delta}{n^2} \\
&\ge & \inf\left\{f + g\right\} + \displaystyle\frac{\delta}{n^2} \\
&> & \left(f + g\right)\left(x_0\right) - \displaystyle\frac{\delta}{n^2} + \displaystyle\frac{\delta}{n^2}.
\end{eqnarray*}
Moreover $h\left(x_0\right) = \delta\theta\left(0\right) = 0, $
then,
$$\inf\left\{f + g + h : \left\|x - x_0\right\| \ge 1/n\right\} > \left(f + g\right)\left(x_0\right) = \left(f + g + h\right)\left(x_0\right). $$
Thus $\left(g + h\right) \in G_n $ and $G_n$ is a dense subset in $Y. $

Therefore the set $\displaystyle\bigcap_{n \ge 1}G_n $ is residual in $Y. $ Following the proof of [15], we can show $f + g $ attains its strong minimum for each $g \in \displaystyle\bigcap_{n \ge 1}G_n. $ \blacksquare

\end{document}